\documentclass[10pt,twoside]{article}
\usepackage{amsfonts}
\usepackage{amssymb,amsmath}
\begin{document}
\begin{center} \textbf{The integer recurrence}
$P(n)=a+P(n-\phi(a))  $ I \\
\textbf{}
%\vspace{12pt}
\\Constantin M. Petridi
\\ cpetridi@hotmail.com\\
\vspace{20pt}
 \textbf{Abstract}\end{center}
\par
\hspace{0.2in}
\begin{tabular}{p{12cm}}
\begin{small}We prove that for a positive integer a the integer
sequence $P(n)$ satisfying for all $n, -\infty<n<\infty,$ the
recurrence $P(n)=a+P(n-\phi(a))$, $\phi(a)$ the Euler function,
generates in increasing order all integers $P(n)$ coprime to $a$.
The finite Fourier expansion of $P(n)$ is given in terms of $a$,
n, and the $\phi(a)$-th roots of unity. Properties of the sequence
are derived.
\end{small}
\end{tabular}
\vspace{25pt}
\par
\textbf{1. Introduction}
\par
For a positive integer $a=\prod_{i=1}^{\omega}\,p_{i}^{e_{i}},\;$
we set $\;R(a)=\prod_{i=1}^{\omega}p_{i}\;$,
$\;Q(a)=\prod_{i=1}^{\omega}\,(p_{i}-1)\;$,
$\;\phi(a)=\frac{a}{R(a)}Q(a)$, which, if $a$ is fixed, we write
$R,\;Q,\;\phi\;$, respectively. Note that $Q|\phi$ and that
$\phi(R(a))=Q(R(a))$.
\par
We show that the integer sequence $P(n)$ defined for all $n$
running from $-\infty$ to $\infty$ , by the inhomogeneous
recurrence of order $\phi(a)$ $$P(n)=a+P(n-\phi(a))$$
 and certain $\phi(a)+1$ initial conditions,
to be specified, gives in increasing order, $P(n) < P(n+1)$, all
integers coprime to $a$.
\par
The minimal integer recurrence satisfied by $P(n)$ is
$$P(n)=R(a)+P(n-Q(a)).$$
Using known facts from the theory of linear integer recurrences,
$P(n)$ is explicitly expressed as a finite Fourier expansion
involving $a$, $n$, and the $\phi(a)$-th roots of unity.
Similarly $P(n)$ can be expressed as a finite Fourier expansion
involving $R(a),\;n$, and the $Q(a)$-th roots of unity.
\par
Properties of the function $P(n)$ are established, such as f.ex.
$\lim_{n\rightarrow\infty}\frac{P(n)}{n}=\frac{R(a)}{Q(a)}.$
\par
The infinite sequence ${P(n)}$, with index $n$ suitably
normalized, can be regarded as the natural extension to
$\mathbb{Z}$ of the "Euler" set of $a$, namely the $\phi(a)$
positive integers smaller than $a$ and coprime to $a$, arranged
in increasing order.
\vspace{25pt}
\par
\textbf{2. The sequence $P(n)$}
\par
We first prove the minimal integer recurrence satisfied by
$P(n)$, from which, afterwards, we deduce the integer recurrence
$P(n)=a+P(n-\phi(a))$.
\par
 \textbf{Theorem 1.} For a fixed
integer $a\geq 1$, let $a_{1}(=1) < a_{2}<\ldots<a_{Q}(=R-1)$
denote the $Q$ positive integers smaller than $R$ and coprime to
$R$. The integer sequence $P(n)$, normalized for $n=1$ to give
$P(1)=R+1$, which satisfies for all $n$ running from $-\infty$ to
$\infty$, the integer recurrence
$$P(n) = R+P(n-Q),$$ and the $Q+1$ initial conditions
$$P(-Q+1)=a_{1}\;\;\;(=1)$$ $$P(-Q+2)=a_{2}$$
$$\hspace{5cm}\ldots\hspace{1cm}\ldots\hspace{1cm}\ldots\hspace{4.3cm}(1)$$
$$P(0)=a_{Q}\;\;\;(=R-1)$$ $$P(1)=a_{Q+1}\;\;\;(=R+1),$$ generates in
monotonically increasing order, $P(n)<P(n+1)$, all integers
coprime to $a$.
\par
\textbf{Proof.} We first show that $P(n)$ is coprime to $a$.
 \par
For $n=0$ this is obvious since $P(0) = R-1$.
\par
For $n>0$, we have by definition $$P(n)=R+P(n-Q)$$ $$P(n-Q) =
R+P(n-2Q)$$
$$\hspace{4cm}\ldots\hspace{1cm}\ldots\hspace{1cm}\ldots\hspace{4cm} $$
$$P(n-(k-1)Q)=R+P(n-kQ),$$ which, when added, give for any integer
$k\geq 0$ $$\hspace{4.2cm}P(n)=kR+P(n-kQ).\hspace{4.5cm}(2)$$
Setting $k=\big[ \frac{n+Q-1}{Q}\big]$, where $[x]$ the greatest
integer $\leq x$, we have
$$\frac{n-1}{Q}<\Big[\frac{n+Q-1}{Q}\Big]\leq \frac{n+Q-1}{Q}$$
$$n-1<\Big[\frac{n+Q-1}{Q}\Big]Q\leq n+Q-1$$
$$-Q+1\leq n-\Big[\frac{n+Q-1}{Q}\Big]Q<1.$$
Hence from (1) we get $$1\leq\;
P\Big(n-\Big[\frac{n+Q-1}{Q}\Big]Q\Big)\;< R+1.$$ Since
$P\big(n-\big[\frac{n+Q-1}{Q}\big]Q\big)$ is coprime to $a$
because of (1), it follows from (2) that $P(n)$ is also coprime to
$a$.
\par
For $n<0$ it can be shown, by a similar argument, that
$$P(-n)=-kR+P(-n+kQ),$$ and by choosing
$k=\big[\frac{n+1}{Q}\big]$ we again derive that $P(-n)$ is
coprime to $a$.
\par
It results that for any $n\in\mathbb{Z}$ and any $k\in
\mathbb{Z}$, we have
$$\hspace{4.4cm}P(n)=kR+P(n-kQ)\hspace{3.5cm}(3)$$
We now show that $P(n)<P(n+1).$ Substituting in (3) $n$ by $n+1$
and subtracting we get
$$P(n+1)-P(n)=P(n+1-kQ)-P(n-kQ).$$
Choosing, as above, for $k$ either $\big[ \frac{n+Q-1}{Q}\big]$ or
$\big[\frac{n+1}{Q}\big]$, we infer from the initial conditions
(1), that $P(n)<P(n+1)$.
\par
From Theorem 1 we obtain
\par
\textbf{Theorem 2.} For a fixed integer $a\geq 1$, let
$a_{1}(=1)<a_{2}<\ldots <a_{\phi}(=a-1)$ denote the $\phi$
positive integers smaller than $a$ and coprime to $a$. The
integer sequence $P(n)$, normalized for $n=1$ to give $P(1)=a+1$,
which satisfies for all $n$ running from $-\infty$ to $\infty$,
the integer recurrence
$$P(n)=a+P(n-\phi),$$ and the $\phi+1$ initial conditions
$$P(-\phi+1)=a_{1}\;\;\;(=1)$$ $$P(-\phi+2)=a_{2}$$
$$\hspace{4.5cm}\ldots\hspace{1cm}\ldots\hspace{1cm}\ldots\hspace{4.5cm}$$
$$P(0)=a_{\phi}\;\;\;(=a-1)$$ $$P(1)=a_{\phi+1}\;\;\;(=a+1),$$
generates in monotonicaly increasing order $P(n)<P(n+1)$, all
integers coprime to $a$.
\par
\textbf{Proof.} Setting in (3) $k=\frac{a}{R}$, we have
$$\hspace{5cm}P(n)=a+P(n-\phi).\hspace{3.5cm}\Box$$
%\vspace{5pt}
\par
\textbf{3. The Fourier expansion of $P(n)$}
\par
The integer recurrence $P(n)=R+P(n-Q)$ is inhomogeneous.
Subtracting $P(n+1)=R+P(n+1-Q)$ we get the homogeneous recurrence
$$P(n+1)-P(n)-P(n+1-Q)+P(n-Q)=0.$$
Its characteristic polynomial is
$$ x^{Q+1}-x^{Q}-x+1=(x-1)^{2}(x-e^{2\pi i\frac{1}{Q}})\ldots (x-e^{2\pi
i\frac{Q-1}{Q}}).$$ Using known facts from the theory of integer
recurrences, we therefore have
$$P(n)=c_{0}n+c_{1}+\sum_{\nu=1}^{Q-1}c_{\nu+1}e^{2\pi
i\frac{\nu}{Q}n}$$
$$\hspace{4cm}=c_{0}n+\sum_{\nu=0}^{Q-1}c_{\nu+1}e^{2\pi
i\frac{\nu}{Q}n}.\hspace{3.8cm}(4)$$ The $Q+1$ coefficients
$c_{\nu}$ can be determined by solving the system of $Q+1$ linear
equations, resulting from following $Q+1$ initial conditions
$$c_{0}(-Q+1)+\sum_{\nu=0}^{Q-1}c_{\nu+1}e^{2\pi
i\frac{\nu}{Q}(-Q+1)}=a_{1}$$
$$c_{0}(-Q+2)+\sum_{\nu=0}^{Q-1}c_{\nu+1}e^{2\pi
i\frac{\nu}{Q}(-Q+2)}=a_{2}$$
$$\hspace{4.5cm}\ldots\hspace{1cm}\ldots\hspace{1cm}\ldots\hspace{4.5cm}$$
$$c_{0}(0)+\sum_{\nu=0}^{Q-1}c_{\nu+1}e^{2\pi
i\frac{\nu}{Q}(0)}=a_{Q}$$
$$c_{0}(1)+\sum_{\nu=0}^{Q-1}c_{\nu+1}e^{2\pi
i\frac{\nu}{Q}}(1)=a_{Q+1}.$$ Another way is to take only the
first $Q$ equations of above system, transfer the terms
$c_{0}(-Q+1), \ldots ,c_{0}(0)$ to the right side, and find the
value of $c_{0}$ afterwards. Accordingly we write
$$\sum_{\nu=0}^{Q-1}c_{\nu+1}e^{2\pi i\frac{\nu}{Q}(-Q+1)}=a_{1}-c_{0}(-Q+1)$$
$$\sum_{\nu=0}^{Q-1}c_{\nu+1}e^{2\pi i\frac{\nu}{Q}(-Q+2)}=a_{2}-c_{0}(-Q+2)$$
$$\hspace{4.5cm}\ldots\hspace{1cm}\ldots\hspace{1cm}\ldots\hspace{4.5cm}$$
$$\sum_{\nu=0}^{Q-1}c_{\nu+1}e^{2\pi i\frac{\nu}{Q}(0)}=a_{Q}-c_{0}(0).$$
Putting (Vandermonde)
$$D_{a}=|e^{2\pi i\frac{k}{Q}\ell}|,\;\;\; 0\leq k\leq Q-1,\;\;\;-Q+1\leq \ell\leq 0$$
we get for the coefficients
$$\hspace{2.7cm}c_{\nu}=\frac{1}{D_{a}}\;\sum_{\mu=1}^{Q}(-1)^{\mu
-1}(a_{\mu}-c_{0}(-Q+\mu))\;D_{\nu,\mu},\hspace{2.9cm}(5)$$ where
$D_{\nu,\mu}$ are the $(Q-1)\times (Q-1)$ minors of $D_{a}$,
obtained by replacing the $\mu$-th column with
$a_{\mu}-c_{0}(-Q+\mu)$.
\par
In order to find the value of $c_{0}$ we substitute in (4) $n$ by
$n+Q$. This gives
$$P(n+Q)=c_{0}(n+Q)+\sum_{\nu=0}^{Q-1}c_{\nu}e^{2\pi i
\frac{\nu}{Q}(n+Q)}$$
$$=c_{0}Q+c_{0}n+\sum_{\nu=0}^{Q-1}c_{\nu}e^{2\pi i
\frac{n}{Q}n}$$
$$=c_{0}Q+P(n).$$
\vspace{50pt}
 But $P(n+Q)=R+P(n)$, as can be seen from (3), if we
substitute $n$ by $n+kQ$ and put $k=1$.
 Hence, $$c_{0}Q=R$$
$$c_{0}=\frac{R}{Q}.$$
Summarizing, we have proved
\par
\textbf{Theorem 3.} $P(n)$, $-\infty<n<\infty$, can be expressed
by the Fourier expansion
$$P(n)=\frac{R}{Q}n+\sum_{\nu=0}^{Q-1}c_{\nu}e^{2\pi i
\frac{\nu}{Q}n},$$ where the coefficients $c_{\nu}$ are given by
(5).
\par
Exactly the same procedure as above, applied to the
recurrence
$$P(n)=a+P(n-\phi)$$
gives
\par
\textbf{Theorem 4.} $P(n)$, $-\infty<n<\infty$, can be expressed
by the Fourier expansion
$$P(n)=\frac{R}{Q}n+\sum_{\nu=0}^{\phi-1}d_{\nu}e^{2\pi i
\frac{\nu}{\phi}n},$$ where the coefficients $d_{\nu}$, depending
on $\phi$, are given by a formula similar to (4).
\par
Dividing by $n$, we get
\par \textbf{Corollary 1.}$$ \lim_{n\rightarrow
\infty}\frac{P(n)}{n}=\frac{R}{Q}.$$ \par\vspace{5pt}
 Following transformation formulas are
an immediate consequence of above Theorems.
\par
\textbf{Corollary 2.} If $R(a)=R(b)$, then for all $n$
$$P(n-Q(a))=P(n-Q(b))$$
$$P(n-\phi(a))=P(n-\phi(b)).$$
\par
\textbf{Note}
\par
In a next communication, we examine $P(n)$ from the angle of
generating functions. To this end, we introduce the GF of the
sequence ${P(n)}$ for a fixed integer $a\geq 1$.
$$G(t)=\sum_{n=1}^{\infty}P(n)t^{n}.$$
The recurrence $P(n)=a+P(n-\phi)$ promptly gives
$$G(t)=a\frac{t}{(t-1)(t^{\phi}-1)}-\frac{1}{t^{\phi}-1}\sum_{\nu=1}^{\phi}a_{\nu}t^{\nu}.$$
The coefficients are then expressed by
$$P(n)=\frac{a}{2\pi i}\int_{c}\frac{t}{(t-1)(t^{\phi}-1)}\,\frac{1}{t^{n+1}}dt
-\sum_{\nu=1}^{\phi}\frac{a_{\nu}}{2\pi
i}\int_{c}\frac{t^{\nu}}{t^{\phi}-1}\,\frac{1}{t^{n+1}}dt.$$
Evaluation of the integrals or expansion in series of the GF
yields expressions of $P(n)$ in terms of $a$, $n$, $\phi(a)$ and
the "Euler" set $\{a_{\nu}\}$. Comparison with the resp. Fourier
expansions results in identities involving the parameters.
\vspace{25pt}
\par
\textbf{Acknowledgment} We wish to thank Peter Krikelis,
University of Athens, Mathematics department, for his help.
\end{document}